\documentclass[11pt]{article}
\usepackage{amsmath,amssymb,a4wide,amsthm,color,graphicx,verbatim,algorithm,caption,enumerate,multirow}
\usepackage{algorithmic}
\usepackage{refcount}
\usepackage{epstopdf}

\newtheorem{theorem}{Theorem}[section]

\newtheorem{remark}[theorem]{Remark}

\newtheorem{ex}[theorem]{Example}

\newtheorem{definition}[theorem]{Definition}

\newcommand{\st} {\ensuremath{|\;}}

\newcommand{\cl}  {{\rm cl  \,}}
\newcommand{\bd}  {{\rm bd \,}}

\newcommand{\Int} {{\rm int \,}}
\newcommand{\conv}  {{\rm conv \,}}
\newcommand{\cone}{{\rm cone\,}}

\newcommand{\R}{\mathbb{R}}

\newcommand{\recc}{{\rm recc \,}}
\DeclareMathOperator*{\argmin}{arg\,min}

\author{\.Irfan Caner Kaya \thanks{Bilkent University, Department of Industrial Engineering, Bilkent, Ankara, 06800 Turkey} \and Firdevs Ulus \thanks{Bilkent University, Department of Industrial Engineering, Bilkent, Ankara, 06800 Turkey, e-mail: firdevs@bilkent.com.tr}}

\title{An iterative vertex enumeration method for objective space based {vector} optimization algorithms}

\begin{document}
\maketitle
	
\begin{abstract} \noindent
An application area of vertex enumeration problem (VEP) is the usage within objective space based linear/convex {vector} optimization algorithms whose aim is to generate (an approximation of) the Pareto frontier. In such algorithms, VEP, which is defined in the objective space, is solved in each iteration and it has a special structure. Namely, the recession cone of the polyhedron to be generated is the {ordering} cone. We {consider and give a detailed description of} a vertex enumeration procedure, which iterates by calling a modified `double description (DD) method' that works for such unbounded polyhedrons. We employ this procedure as a function of an existing objective space based {vector} optimization algorithm (Algorithm 1); and test the performance of it for randomly generated linear multiobjective optimization problems. We compare the efficiency of this procedure with another existing DD method as well as with the current vertex enumeration subroutine of Algorithm 1. We observe that the modified procedure excels the others especially as the dimension of the vertex enumeration problem (the number of objectives of the corresponding multiobjective problem) increases.
		
\medskip
		
		
\medskip
		
		
\end{abstract}
	
\section{Introduction}
\label{sect:intro}
A polyhedron $P \subseteq \R^d$ can be represented as intersection of finitely many halfspaces or as convex hull of its vertices added to the conic hull of its extreme directions. The problem of computing the vertex representation of $P$ from its halfspace representation is called the \textit{vertex enumeration problem} (VEP). VEP has been studied for many years, starting at latest from the 1950s, see for instance~\cite{dd,bal61}. The difficulty of the problem is known~(\cite{davis3,kac08}) and there are many studies to propose efficient algorithms or to improve the efficiency of existing ones (\cite{matthschmidt80,dyerproll1982,fukpro,avisonly98}).

The vertex enumeration problem as defined above is sometimes called the \emph{off-line} VEP, whereas finding the vertices of a polyhedron $P^{''}$ given as intersection of a single halfspace $H$ and another polyhedron $P^{'}$ whose vertex representation is known is called the \emph{on-line} VEP. Note that an on-line vertex enumeration algorithm can be called repetitively in order to solve an off-line VEP. Different from those, there are also (simplex-type) pivoting algorithms that solve off-line VEP directly, see for instance~\cite{bal61,alt75,avisfukuda92}.

The problem of finding vertices of a polyhedron is fundamental and it is a base of some other algorithms including the simplex algorithm to solve linear programs and outer approximation algorithms for DC programming problems, see for instance~\cite{dc_tuy}. It has also been an important part of objective space based algorithms designed to solve linear (\cite{ben1,lvop,csirmaz16}) or convex (\cite{ehrshao,cvop}) multiobjective {and vector optimization problems (VOPs)}. In general, these algorithms aim to find (a polyhedral approximation of) the set of all nondominated points in the objective space, namely the Pareto frontier. Clearly, for a VOP with $d$ objectives, the objective space is $\R^d$. In each iteration of such algorithm, an on-line vertex enumeration problem of size $d$ has to be solved. Note that for a {VOP with an ordering cone $K$}, the nondominated {(K-minimal)} points of the image of the feasible region is the same as {that} of the \emph{upper (extended) image}, which is {$K$} added to the image of the feasible region. Indeed, for these algorithms, the idea is to approximate the upper image by inner and outer polyhedral sets. 

{Solving sequential online VEPs is also used in order to solve data envelopment analysis (DEA) problems recently. Recall that DEA is used to compute the efficiency of decision-making units (DMUs) with common inputs and outputs (for a review, see \cite{cooper2011handbook,cook2009data}), and it has many application areas including, banking, healthcare, energy and environmental sciences, agriculture, see for instance the survey papers \cite{zhou2008survey,liu2013survey}. In \cite{ehrgott2019multiobjective}, Ehrgott, Hasannasab and Raith  proposed an algorithm in order to generate the extreme points and facets of the efficient frontier of data envelopment analysis (DEA) problems using the geometric duality theory for linear multiobjective optimization problems (MOPs). Accordingly, instead of relying on solving a linear program for each DMU as many other DEA solution approaches from the literature do, their algorithm is based on solving online VEPs in each iteration and does not solve any LP. It is demonstrated in~\cite{ehrgott2019multiobjective} that their algorithm is computationally comparable with the standard DEA approach for some real life problems and it is faster than that for large-scale artificial data sets for which the percentage of efficient DMUs is rather small.}

Non-pivoting on-line vertex enumeration algorithms in the literature are generally designed to find the vertices of bounded polyhedrons, see for instance~{\cite{chen90}}. However, for {aforementioned vector} optimization algorithms, one needs to find the vertices of unbounded polyhedrons since the upper image of a VOP is an unbounded polyhedron whose recession cone is (at least) the ordering cone. {In 2010, using projective geometry, Burton and Ozlen \cite{burton2010projective} proposed a vertex enumeration method which works for unbounded polyhedrons with known recession cones.} %In this study, we revisit this vertex enumeration procedure from~\cite{burton2010projective} and provide a detailed description of it. 

{In this study, we first consider the `standard' double description (DD) method proposed for bounded polyhedrons. The detailed description is provided in~\cite{fukpro}, where the problem is given in an equivalent form of enumerating the extreme directions of a polyhedral cone.} Note that the DD method is designed to solve the on-line vertex enumeration problem and by employing it in each iteration, one can solve an off-line VEP as well. In particular, one needs to start with a sufficiently large bounded polyhedron $P^0$ that contains the set of all vertices of $P$. This method can be applied both for bounded or unbounded $P$ for which the recession cone is not necessarily known. As long as the vertices of $P$ are known to be included in $P^0$, one can use DD method iteratively and at the end of the final iteration, one needs to get rid of the vertices that are on the boundary of $P^0$. The main difficulty in this approach is to find such $P^0$. Also, taking $P^0$ too large may result in numerical issues when implemented.

{We then consider the modified DD method (proposed in \cite{burton2010projective}), which works for unbounded $P$ whose recession cone, say $K$, is known. We provide a detailed description of the algorithm in the Euclidean space (instead of the oriented projective space). Note that} the main difference of the modified method is the use of extreme directions of {the current approximation} $P^{'}$ in the computations of the vertices of the updated polyhedron $P^{''}$. Clearly, the modified DD method can also be called iteratively in order to solve an off-line VEP. This time one needs to start with a single vertex $v_0 \in \R^d$ such that the initial polyhedron {$P^0:=v^0+K$} contains the polyhedron to be computed. {For computing the upper image of a multiobjective optimization problem (MOP), where the ordering cone is the nonegative cone,} \emph{ideal point} is the natural candidate for this initial vertex $v^0$.

We implement the {standard and the modified DD methods, and compared their efficiencies through computational tests. Indeed,} we employ both methods as a subroutine within a MATLAB implementation of the objective space based convex {vector} optimization algorithm proposed in~\cite{cvop}. The current implementation of the algorithm calls in each iteration a vertex enumeration procedure (\emph{vert}) which mainly uses the `qhull' function of MATLAB (\cite{qhull}). Note that this vertex enumeration procedure is first employed within the MATLAB implementation of \emph{bensolve}, which is a linear vector optimization solver, see~\cite{bensolve12,bensolve2}. We test the performances of the original and the modified DD methods together with \emph{vert} through randomly generated linear MOPs. 

In Section~\ref{sect:Prelim}, we provide preliminaries on basic convex analysis and on convex {VOPs} as well as the convex VOP algorithm proposed in~\cite{cvop}. The DD method and the modified variant are described in Section~\ref{sec:DD}. The variants of the convex VOP algorithm using different DD methods are explained in Section~\ref{sec:alg12}. The computational results are presented and discussed in Section~\ref{sec:numerical}. We conclude our discussion in Section~\ref{sec:concl}.

\section{Preliminaries}
\label{sect:Prelim}

The boundary, the interior, the convex hull and the conic hull of a subset $S \subseteq \R^d$ is denoted by  $\bd S$, $\Int S$, $\conv S$, $\cone S$, respectively.

Let $S$ be a convex subset of $\R^d$ and $F\subseteq S$ be a convex subset. If $\lambda y^1+(1-\lambda) y^2\in F$ for some $0 < \lambda <1$ holds only if both $y^1$ and $y^2$ are elements of $F$, then $F$ is a \textit{face} of $S$. A zero dimensional face is an \emph{extreme point}, a one dimensional face is an \emph{edge} and a $d-1$ dimensional face is a \emph{facet} of $S$, see~\cite{rockafellar}.

For a subset $S$ of $\R^d$, $z\in\mathbb{R}^{d}\setminus \{0\}$ is a \textit{recession direction} (or simply \emph{direction}) of $S$, if $y+ \gamma z\in S$ for all $\gamma \geq 0$ and $y\in S$. The set of all recession directions constitute the \textit{recession cone} of $S$ which is denoted by $\recc S$. 
{A recession direction $z \in \recc S \setminus \{0\}$ of convex set $S$ is said to be an \textit{extreme direction} of $S$, if $\{v+rz \in \R^d \st r \geq 0\}$ is a face for some extreme point $v$ of $S$. $S \subseteq \R^d$ is \textit{bounded} if $\recc S = \{0\}$.}

{A closed convex pointed solid cone $K\subseteq \R^d$ defines a partial order on $\R^d$ as follows: $a\leq_K b \iff b-a \in K$.} Throughout, $\R^d_+:=\{y\in \R^d \st y_i \geq 0, i=1,\ldots,d\}$ is the nonnegative cone in $\R^d$, $e^i$ is the unit vector in $\R^d$ with $i^{th}$ component being $1$.%, and $e\in\R^d$ is the vector of ones.

\subsection{Representations of a Convex Polyhedron}\label{subsec:prelim_1}

If a convex set $P$ can be written as $P = \{y\in \mathbb{R}^{d}\ | \ A^Ty\geq b\}$,
where $A\in \mathbb{R}^{d\times k}$ and $b\in \mathbb{R}^{k}$, then it is called a \textit{polyhedral convex set} or a \textit{convex polyhedron}. Note that $P$ is intersection of finitely many half-spaces, namely,
\begin{equation} \label{Hrep:A}
P =\bigcap\limits^k_{i=1} \{ y\in \mathbb{R}^{d}\ | \  a_i^Ty\geq b_i\},
\end{equation}
where $a_i \in \R^d$ is the $i^{th}$ column of matrix $A$ and $b_i\in\R$ is the $i^{th}$ component of $b$. The representation given in~\eqref{Hrep:A} (with the assumption that there are no redundant inequalities) is called \textit{H-representation} or \textit{halfspace representation} of $P$. On the other hand, if $P$ has at least one extreme point, it can also be represented as the convex hull of all its extreme points added to the conic hull of all its extreme directions. To be more precise, let $V$ be the finite set of all extreme points (vertices) of $P$ and $D$ be the finite set of  all extreme directions of $P$. Then, $P$ can be written as
\begin{equation} \label{Vrep:A}
P=\conv V + \cone \conv D. 
\end{equation}
The representation given by~\eqref{Vrep:A} of $P$ is called the \textit{V-representation} or the \textit{vertex representation} of the polyhedral convex set $P$. The problem of finding the \textit{V-representation} of a set given its \textit{H-representation} is called the \textit{vertex enumeration problem}.

\subsection{Convex {Vector} Optimization and an Approximation Algorithm} \label{subsec:cvop}

A convex {vector} optimization problem is given by 
\begin{align*} \label{(P)}
\text{minimize~} F(x), {\text{~with respect to~} \leq_K,} \text{~subject to~} x \in \mathcal{X}, 
\tag{P}
\end{align*}
where {$K\in\R^d$ is solid, pointed, polyhedral convex ordering cone, $F:\R^n\to \R^d$ is a $K$-convex continuous function} %$F(x) = (f_1(x), \ldots, f_d(x))^T$ with $f_i: \R^n \rightarrow \R$ for all $i=1,\ldots,d$ are convex functions 
and the feasible set $\mathcal{X} \subseteq \R^n$ is a convex set. The image of the feasible set is given by $F({\mathcal{X}}) = \{ F(x) \in \R^d \ | \ x \in {\mathcal{X}}\}$ and the set {$\mathcal{P} := \cl(F(\mathcal{X})+K)$} is called the \emph{upper (extended) image} of~\eqref{(P)}. It is known that $\mathcal{P}$ is convex and closed. {\eqref{(P)} is said to be a bounded problem if there exists $y\in \R^d$ such that $\mathcal{P}\subseteq \{y\}+K$.}

{If the ordering cone is the nonnegative cone, then \eqref{(P)} is a multiobjective optimization problem and} the \emph{ideal point} of problem~\eqref{(P)} can be found by minimizing {each objective $F_i$}, for $i=1,\ldots,d$ over feasible set $\mathcal{X}$ as long as the corresponding single objective optimization problems have finite optimal objective values. More specifically, let $y_i:=\min\{F_i(x)\ | \ x\in \mathcal{X}\}$. Then, $y^I:=(y_1,\ldots,y_d)^T$ is the ideal point of~\eqref{(P)}. {The problem is bounded if $y^I \in \R^d$.}

For {VOPs}, there are different solution concepts as there is not necessarily a unique `solution' that minimizes all the objective functions simultaneously. Some of the solution concepts are as follows: A point $y \in F(\mathcal{X})$ in the image set is said to be a \emph{{$K$-minimal} (non-dominated)} point if {$$(\{y\}-K) \cap F(\mathcal{X}) = \{y\}.$$} Similarly $y \in F(\mathcal{X})$ is said to be a \emph{{weakly $K$-minimal} (weakly non-dominated)} point if {$$(\{y\}-\Int K)\cap F(\mathcal{X}) = \emptyset.$$} A feasible point $x \in \mathcal{X}$ is said to be a \emph{(weakly) efficient solution} if $F(x)$ is a (weakly) non-dominated point of $F(\mathcal{X})$.

In some applications of {VOPs}, it is important to generate the set of all (weakly) non-dominated points of $F(\mathcal{X})$, which is a subset of the boundary of the upper image. When the problem is linear, then it is possible to generate (the set of all extreme points of) the upper image, see for instance~\cite{ben1,lvop}. If the problem is nonlinear convex, it is not possible to generate the set of all (weakly) non-dominated points in general. Instead, there are objective space based algorithms that can generate polyhedral approximations to the upper image as in~\cite{ehrshao,cvop}. 

The general idea of such an algorithm is as follows. It starts with finding {an initial outer approximation $P^0$ of $\mathcal{P}$. In case of a MOP, this can be done by finding} {the ideal point} $y^I$ of problem~\eqref{(P)}. Then, the initial outer approximation of $\mathcal{P}$ is $P^0 := \{y^I\} + \R^d_+$. 

At $i^{th}$ iteration of the algorithm, the first step is to find the vertices of the current outer approximation $P^{i-1}$. Next, for a vertex $v$ of $P^{i-1}$, single objective convex optimization problem, namely the Pascoletti-Serafini scalarization (\cite{PS}), given by
\begin{align} \label{Pv}
\text{minimize~~~} \alpha \text{~~~subject to~~~} F(x){\leq_K} v+\alpha {k}, \:\: x \in \mathcal{X}, \tag{P($v$)}
\end{align}
{for some fixed $k\in \Int K$,} is solved. Note that this problem finds point $y^v:=v + \alpha^v {k}$ on $\bd \mathcal{P}$ that is ``closest" (through the fixed direction {$k$}) to $v$, where $\alpha^v$ is the optimal objective function value of the program. Moreover, optimal solution $x^v$ is known to be weakly efficient.

If $\alpha^v>\epsilon$, where $\epsilon$ is the predetermined approximation error, then by using the dual solution of this scalar convex optimization problem, one finds a supporting hyperplane of $\mathcal{P}$ at $y^v$. More specifically, if $w \in \R^d$ denotes the dual variable corresponding to the first set of constraints {$F(x)\leq_K v+\alpha k$}; and $w^v$ is the dual optimal solution, then
$$h^v:=\{y\in \R^d \st (w^v)^Ty = (w^v)^Ty^v \}$$ 
supports $\mathcal{P}$ at $y^v$ and $$H^v := \{y\in \R^d \st (w^v)^Ty \geq (w^v)^Ty^v\}$$ is the corresponding halfspace that  contains $\mathcal{P}$ (\cite{cvop}). After computing $H^v$, the outer approximation of the upper image is updated as $P^i := P^{i-1} \cap H^v$ and the $i^{th}$ iteration is completed. 

If $\alpha^v \leq \epsilon,$ then the algorithm continues checking other vertices of the current outer approximation. The algorithm stops when all the vertices are in $\epsilon$-distance to the upper image. One can see the books~\cite{ehrgott_book} and~\cite{lohne} for details of the multiobjective/vector optimization theory and of objective space based (also referred as Benson-type (\cite{ben1})) algorithms. We provide the pseudo-code of the algorithm as proposed in~\cite{cvop}, see {Algorithm 1}.
\begin{algorithm}[!]
	\begin{algorithmic}[1]
		\STATE {Compute an initial outer approximation $P^0$ of $\mathcal{P}$, initialize $\bar{\mathcal{X}}$ accordingly;} \STATE $i=0$;
		\REPEAT
		\STATE continue = 0;
		\textbf{	\STATE~} \hspace{-0.3cm}Compute the vertices $V^i$ of $P^i$;
		\FORALL{$v\in V^i$}
		\STATE Solve \eqref{Pv}, let optimal solution be $(x^v,\alpha^v)$ and dual optimal solution be $w^v$;
		\STATE $\bar{\mathcal{X}} \gets \bar{\mathcal{X}} \cup \{x^v\}$
		\IF{$\alpha^v > \epsilon$}
		\STATE $i\gets i+1$;
		\STATE $P^i = P^{i-1} \cap H^v$;
		\STATE continue = 1;
		\STATE break;
		\ENDIF
		\ENDFOR
		\UNTIL{continue = 0}
		\textbf{\RETURN~} \hspace{-0.3cm} $\bar{\mathcal{X}}$: Set of weakly efficient solutions \\ ~~~~~~~~~~~~$P^i$: Polyhedral outer approximation to the upper image.\hrule
	\end{algorithmic} 
	\caption{An objective space based convex {VOP} Algorithm}
	\label{alg:benson} 
	{ For a MOP, the initialization can be performed as follows:
		{\small  \begin{enumerate}[1:]
				\item Compute $y_i:=\min\{f_i(x)\ | \ x\in \mathcal{X}\}$, $x^i:=\argmin\{f_i(x)\ | \ x\in \mathcal{X}\}$;\vspace{-0.15cm}
				\item $y^I:=(y_1,\ldots,y_d)^T$, $P^0 = \{y^I\} + \R^d_+$, $\bar{\mathcal{X}}=\{x^1, \ldots, x^d\}$;
		\end{enumerate}}
	}
\end{algorithm}

There are different variants of the algorithm in the literature. For example, in each iteration of the variant proposed in~\cite{ehrshao}, the direction parameter ({$k$} in \eqref{Pv}) is chosen in a different way depending on $v$, and the supporting hyperplane is constructed using the differentials of the objective functions instead of using the dual optimal solution of \eqref{Pv}. 

\begin{remark}\label{rem:Benson}
	Note that at each iteration of the algorithm, the first step is to solve an on-line VEP and these VEPs are in a special form. The polyhedron to be found is unbounded since $\mathcal{P}$ is an unbounded set. Moreover, assuming that the {VOP is bounded}, the recession cone of the polyhedron is the {ordering} cone (and not larger than that). Hence the extreme directions of the recession cone are {the extreme directions of $K$, namely $Z:=\{z^1,\ldots,z^m\}$. In case of a MOP,} they are nothing but the unit directions, $\{e^1,\ldots,e^d\}$.
\end{remark}

\section{The Double Description (DD) Method} \label{sec:DD}

We first describe the DD method which works for bounded polyhedrons and then describe a modification of it which works for unbounded polyhedrons with {known recession cones}.

For both methods, let $P^{'}$ be a convex polyhedron, $H$ be a halfspace given by $H:=\{x\in\R^{d}\:|\: a^Tx\geq b\}$ for some $a\in\R^{d}\setminus\{0\}$ and $b\in\R$, and $h:=\bd H$ be the hyperplane given by the boundary of $H$.

\subsection{The DD Method} \label{subsect:DD}	
In this section, we describe the DD method, as provided in~{\cite{burton2010projective,csirmaz16}}.	

Let $P^{'}$ be bounded, $V^{'}$ be the set of its vertices and $F^{'}$ be the set of its facets. The following is a useful definition in order to describe the method.

\begin{definition} \label{defn:adj_v_f}
	If vertex $v\in V^{'}$ is on facet $f\in F^{'}$; then $v$ is said to be an \emph{adjacent vertex} of facet $f$ and $f$ is said to be an \emph{adjacent facet} of vertex $v$.
\end{definition}
For a vertex $v$, let $F_v$ denote the set of all adjacent facets of $v$, and for a facet $f$, let $V_f$ denote the set of all adjacent vertices of $f$. These sets are called the adjacency lists. It is assumed that the sets $V^{'}, F^{'}$ as well as $V_f, F_v$ for all $v\in V^{'}$ and $f \in F^{'}$ are known. 

As it is an important subroutine in DD method, we first describe a procedure to check if there is an edge between given two vertices of polyhedron $P^{'}$. Let $v^{+}, v^{-} \in V^{'}$ be two vertices. In order to check if there is an edge between them, one considers the set of facets which are both adjacent to $v^{+}$ and $v^{-}$. Then, for each facet in this set, one considers the adjacent vertices of it. If the intersection of the set of vertices over all these facets consists of only $v^{+}$ and $v^{-}$ then, there is an edge between the two; otherwise, there is no edge between them, {see \cite[Lemma 5.3]{burton2010projective} and \cite{fukpro}}. Procedure~\ref{alg:isedge} is the pseudo-code for the $isedge$ function, which takes polyhedron $P^{'}$ and vertices $v^{+}$, $v^{-}$ as its input; and returns the set of facets that contains the edge between them if there is any or returns empty set otherwise. Note that with $P^{'}$ being an input we mean that there is an access to $V^{'}, F^{'}$ as well as $V_f$ and $F_v$ for all $f\in F^{'}$ and $v \in V^{'}$. This will be the case for all the procedures described later as well. 

\floatname{algorithm}{Procedure}
\begin{algorithm}[!] 
	\addtocounter{algorithm}{-1}
	\begin{algorithmic}[1]  
		\STATE Let $F^{\pm}:=F_{v^{+}} \cap\ F_{v^{-}}$;
		\STATE Let $V^{\pm}:= V^{'}$;
		\FOR{$i=1:\left| F^{\pm} \right|$}
		\STATE Let $f^i$ be the $i^{th}$ facet in $F^{\pm}$;
		\STATE $V^{\pm} \gets V^{\pm}\cap V_{f^i}$;
		\ENDFOR
		\IF{$V^{\pm} = \{v^{+},v^{-}\}$}
		\RETURN $F^{\pm}$
		\ELSE
		\RETURN $\emptyset$
		\ENDIF
	\end{algorithmic}
	\caption{$isedge(P^{'},v^{+},v^{-})$}
	\label{alg:isedge}
\end{algorithm}

The idea of the double description method is as follows: First, it checks if each vertex $v$ in $V^{'}$ is in the interior of $H$, on the boundary $h$ of $H$, or not included in $H$. Clearly, the ones that are not in $H$ will not be a vertex of the updated polyhedron anymore. As long as there exists at least one vertex that is included in $H$ and there exists at least one vertex that is not included in it then, the algorithm considers each couple of vertices $v^{+}$ and $v^{-}$ in $V^{'}$, where $v^{+} \in \Int H$ and $v^{-} \notin H$, and checks if there is an edge between $v^{+}$ and $v^{-}$. For the couples that form an edge, a new vertex is found by intersecting the edge with hyperplane $h$. This new vertex is a vertex of the updated polyhedron $P^{''}$ and $h$ is a facet of $P^{''}$. Each time a new vertex is found, the adjacency lists are updated accordingly. 

The pseudo-code for the double description method is given by Procedure~\ref{alg:onlinevert}. The function $onlinevert$ takes polyhedron $P^{'}$ and halfspace $H$ as its input and returns the updated polyhedron $P^{''}= P^{'}\cap H$.

\begin{algorithm}[!]
	\begin{algorithmic}[1]
		\STATE Initialize $V^0,V^{+},V^{-}:= \emptyset$; 
		\FORALL{$v\in V^{'}$}
		\IF{$a^Tv > b$ (i.e. $v \in \Int H$)}
		\STATE $V^{+} \gets V^{+} \cup \{v\}$;
		\ELSIF{$a^Tv = b$ (i.e. $v \in h$)} 
		\STATE $V^{0} \gets V^{0} \cup \{v\}$;
		\ELSIF{$a^Tv < b$ (i.e. $v \notin H$)} 
		\STATE $V^{-} \gets V^{-} \cup \{v\}$;
		\ENDIF 
		\ENDFOR
		\IF{$V^{+}\cup V^0 = V^{'}$}
		\RETURN $P^{''} = P^{'}$;
		\ELSIF{$V^{-}  = V^{'}$}
		\RETURN $P^{''} = \emptyset$;
		\ELSE
		\STATE $F^{'} \gets F^{'}\cup \{h\}$, $V_h = \emptyset$;
		\FORALL{$v\in V^{0}$}
		\STATE $V_h \gets V_h \cup \{v\}$ and $F_v \gets F_v \cup \{h\}$;
		\ENDFOR \label{alg:forloopstarts}
		\textbf{\FORALL{$v^{+}\in V^{+}$} 
			\FORALL{$v^{-}\in V^{-}$}
			\IF{$F^{\pm}:= isedge(P^{'},v^{+},v^{-}) \neq \emptyset$}
			\STATE \textnormal{Find} $v^{'}:=[v^{+},v^{-}]\cap h$ 
			\IF{$v^{'}\notin V^{'}$}
			\STATE $V^{'}\gets V^{'}\cup \{v^{'}\}$, $V_h \gets V_h \cup \{v^{'}\}$ \textnormal{and} $F_{v^{'}} = F^{\pm}\cup\{h\}$;
			\ELSE
			\STATE $F_{v^{'}} \gets F_{v^{'}} \cup F^{\pm}\cup \{h\}$ \textnormal{and} $V_h \gets V_h \cup \{v^{'}\}$;
			\ENDIF
			\FORALL{$f \in F_{\pm}$}
			\STATE $V_f \gets V_f \cup \{v^{'}\}$
			\ENDFOR
			\ENDIF
			\ENDFOR
			\ENDFOR}
		\ENDIF	
		\STATE $V^{''} = V^{'}\setminus V^{-}$ and $F^{''} = \cup_{v \in V^{''}}F_v$;
		\STATE $V_f \gets V_f \setminus V^{-}$ for $f \in F^{''}$;
		\RETURN $P^{''}$ with vertices $V^{''}$, facets $F^{''}$ and respective adjacency lists.
		%}
	\end{algorithmic}
	\caption{$onlinevert(P^{'},H)$}
	\label{alg:onlinevert} 
\end{algorithm}

\subsection{The Modified DD Method}\label{subsect:modifiedDD}

We {consider} a modified  double description method which works for unbounded polyhedrons {with known recession cones. This is proposed and described in \cite{burton2010projective}, using oriented projective geometry. Here, we provide a detailed description in Euclidean space.} Let $V^{'}$ be the set of vertices, $F^{'}$ be the set of facets and {$Z = \{z^1,\ldots,z^m\}$} be the set of extreme directions of $P^{'}$. We assume that the recession cone of the updated polyhedron $P^{''}=P^{'}\cap H$ is also {$\cone \conv Z =: K$}. Note that this is the case if this method is used to compute the vertices of (an approximation of) the upper image of a VOP, see Remark~\ref{rem:Benson}. 

In order to describe the modified DD method, in addition to Definition~\ref{defn:adj_v_f}, we need the following:
\begin{definition}
	Let $P$ be an unbounded polyhedron. An extreme direction $z\neq 0$ of $\recc P$ is said to be an \emph{adjacent direction} of facet $f$ of $P$ if there exists an extreme point $v$ of $P$ such that $\{v+\gamma z \st \gamma \geq 0\}$ forms an edge of $P$ and is on facet $f$. Symmetrically, $f$ is said to be an \emph{adjacent facet} of direction $z$.
\end{definition}

\begin{remark} \label{rem:notation_Vf}
	As before, $F_v$ denotes the set of all adjacent facets of vertex $v$. Different from the previous case, for this method, $V_f$ denotes the set of adjacent vertices together with the set of adjacent directions of facet $f$. Moreover, $F_z$ is the set of adjacent facets of extreme direction $z$. {Here, the extreme directions are treated as vertices. Indeed, they are vertices of the polyhedron when it is described using the oriented projective geometry, see \cite{burton2010projective}.} Hence, from now on whenever we mention adjacent vertices of a facet, we mean the union of adjacent vertices and adjacent directions of it.
\end{remark}

Note that $isedge$ function for given two vertices of polyhedron $P^{'}$ does not use the coordinates of the vertices but only the adjacency lists. Then, by definition of an adjacent facet of an extreme direction $z$ and by the usage of notation $V_f$ (the union of adjacent vertices and adjacent directions of facet $f$), $isedge(P^{'},z,v)$ returns the set of facets on which $\{v+\gamma z \st \gamma \geq 0\}$ lays if this is an edge of $P^{'}$; and returns empty set otherwise, {see \cite[Theorem 5.6]{burton2010projective} for the proof.}

When treating the extreme directions as vertices, one needs to be careful whenever hyperplane $h$ is parallel to some of these extreme directions. Note that if $v^{-}$ is not in $H$, $h$ is not parallel to $z$ and $\{v^{-}+\gamma z \st \gamma \geq 0\}$ is an edge of $P^{'}$; then, $h$ intersects with this edge at a singleton, namely at $v^{'}:= \{v^{-}+\gamma z \st \gamma \geq 0\}\cap h$. Clearly, $v^{'}$ is a vertex of the updated polyhedron. 

If $h$ is parallel to a direction $z$, it does not intersect any edge of the form $\{v+\gamma z \st \gamma \geq 0\}$. Instead, (as long as it cuts) it cuts polyhedron $P^{'}$ in parallel to direction $z$. Hence, $h$ is an adjacent facet of direction $z$. Indeed, there must be a vertex $v$ of $P^{''}$ such that $\{v+\gamma z \st \gamma \geq 0\}$ is an edge of $P^{''}$ and lays on $h$. Similarly, $z$ is an adjacent direction of facet $h$.  

The general structure of the modified DD method ($onlinevert2$) is similar to $onlinevert$. The lines between 1-19 and 36-38 of Procedure~\ref{alg:onlinevert} are exactly the same for $onlinevert2$ as well. The only difference is in the main loop and its pseudo-code can be seen in Procedure~\ref{alg:onlinevert2}. 

\floatname{algorithm}{Procedure}
\begin{algorithm}[h]
	\begin{algorithmic}[1]
		%	{\footnotesize 
		\setcounterref{ALC@line}{alg:forloopstarts}
		\FORALL{$v^{+}\in V^{+}\cup Z$}
		\IF{$v^{+}\in Z$ and $a^Tv^{+} = 0$ ($h$ is parallel to $v^{+}$)}
		\STATE $V_h \gets V_h \cup \{v^{+}\}, F_{v^{+}} = F_{v^{+}} \cup \{h\}$;
		\ELSE
		\FORALL{$v^{-}\in V^{-}$}
		\IF{$F^{\pm}:= isedge(P^{'},v^{+},v^{-}) \neq \emptyset$}
		\IF{$v^{+}\in V^{+}$}
		\STATE Find $v^{'}:=[v^{+},v^{-}]\cap h$;
		\ELSE
		\STATE Find $v^{'}:=\{v^{-}+\gamma v^{+} \st \gamma \geq 0\}\cap h$;	
		\ENDIF
		\IF{$v^{'}\notin V^{'}$}
		\STATE $V^{'}\gets V^{'}\cup \{v^{'}\}$, $V_h \gets V_h \cup \{v^{'}\}$ and $F_{v^{'}} = F^{\pm}\cup\{h\}$;
		\ELSE
		\STATE $F_{v^{'}} \gets F_{v^{'}} \cup F^{\pm}\cup \{h\}$ and $V_h \gets V_h \cup \{v^{'}\}$;
		\ENDIF
		\FORALL{$f \in F_{\pm}$}
		\STATE $V_f \gets V_f \cup \{v^{'}\}$;
		\ENDFOR
		\ENDIF
		\ENDFOR
		\ENDIF
		\ENDFOR	
	\end{algorithmic}
	\caption{$onlinevert2(P^{'},H)$ (substitution to lines 20-34 of Procedure~\ref{alg:onlinevert})}
	\label{alg:onlinevert2} 
\end{algorithm}
The main loop goes over all vertices that are in $\Int H$ and over all directions $Z$ (line 20). If $v^+$ is one of the directions, say $z \in Z$, and if $z$ is parallel to $h$, then $z$ is added as an adjacent `direction' (vertex) of $h$ and $h$ is added as an adjacent facet of $z$ (lines 21-22). Otherwise, i.e., when $v^+ \in V^{+}$ or when $v^{+}$ is a direction which is not parallel to $h$; the algorithm goes over all vertices $v^{-}$ that are not included in $H$ and checks if there is an edge formed by $v^{+}$ and $v^{-}$ (line 25). If there exists an edge and if $v^{+}$ is a vertex (not a direction), then the line segment $[v^+,v^-]$ intersects $h$ at a single point $v^{'}$ (line 27). If $v^{+}$ is a direction and together with $v^{-}$ it forms an edge, then $\{v^{-}+\gamma v^{+} \st \gamma \geq 0\}$ intersects $h$ at a single point $v^{'}$ (line 29). In both cases, $v^{'}$ is a vertex of the updated polyhedron. The adjacency lists are updated in the same way as it is done for $onlinevert$ (lines 31-38).
As a final step after the main loop, the final set of vertices and facets together with adjacency lists are updated as it is done in $onlinevert$ (see lines 36-37 of Procedure~\ref{alg:onlinevert}). Then, the updated polyhedron $P^{''}$ is returned (line 38 of Procedure~\ref{alg:onlinevert}). 

\section{{VOP} Algorithms with DD Methods} \label{sec:alg12}

In the current implementation of Algorithm 1~(\cite{cvop}), vertex enumeration problems are solved using a MATLAB function (\emph{vert}) written for MATLAB implementation of an objective space based (Benson-type) linear {vector} optimization solver \emph{bensolve} (\cite{bensolve12}). Even though an on-line vertex enumeration problem is solved at each iteration of Benson-type algorithms, \emph{vert} solves an off-line vertex enumeration problem. Hence, at each iteration, it computes all the vertices from an H-representation of the current outer approximation even though many vertices are already found in earlier iterations.

The double description methods described in Section~\ref{sec:DD} are used in order to solve the on-line vertex enumeration problem. For Algorithm~\ref{alg:benson}, online vetex enumeration subroutine can be called repetitively in order to compute the vertices of the outer approximation in each iteration. Below, we describe two variants of Algorithm~\ref{alg:benson} that are using $onlinevert$ and $onlinevert2$, respectively.

{Throughout, we assume that the double description of the ordering cone $K$ is known. In other words, adjacent facets $\tilde{F}_{z^i}$ for each extreme direction $z^i \in Z=\{z^1,\ldots,z^m\}$ and adjacent directions $\tilde{V}_{\tilde{f}^{i}}$ for each facet, say $\tilde{f}_i \in \tilde{F}:=\{\tilde{f}_1,\ldots,\tilde{f}_l\}$ are known. Note that for $K=\R^d_+$, we have $Z=\{e^1,\ldots,e^d\}$ and conic hull of any $d-1$ extreme directions from $Z$ forms a facet of $\R^d_+$. More specifically, the set of all facets are $\tilde{F}:=\{\tilde{f}_1,\ldots, \tilde{f}_d\}$ with $\tilde{f}_i := \cone \conv \{Z\setminus\{e^i\}\}$. Then, the adjacent directions of facet $\tilde{f}_i$ is $\tilde{V}_{\tilde{f}^i} = Z\setminus \{e^i\}$ and the adjacent facets of direction $e^i$ is $\tilde{F}_{e^i} = \tilde{F}\setminus\{\tilde{f}^i\}$. For the computational tests in Section~\ref{sec:numerical}, we consider MOPs where $K=\R^d_+$.}

\subsection{Variant 1: {VOP} Algorithm with DD Method} \label{subsec:alg1}

The double description method is used when the initial polyhedron is bounded. However, it can still be used in order to compute the vertices of unbounded polyhedrons. In order to do that, one needs to initialize the algorithm with a large enough initial polyhedron which guarantee to include all the vertices of the polyhedron. Note that as there are finitely many vertices, there exists a bounded polyhedron which contains all. In general any polyhedron satisfying this property can be taken as the initial one. 

Recall that Algortihm~\ref{alg:benson} starts by finding {an initial outer approximation $P^0$. For a bounded problem this can be described as $\{y\}+K$ for some $y\in \R^d$ and for a MOP, $y$ can be taken as the ideal point $y^I$.} Moreover, we assume that there exists a sufficiently large number $M$ such that the set of all nondominated points is a subset of {$P^0:=\{y\}+\conv\{0, Mz^1, \ldots, Mz^m\}$}. Indeed, if the feasible region of the problem is compact, this would be the case as the image of the feasible region in the objective space is bounded. For linear {MOPs}, there exists such $M$ as long as the ideal point is finite (which may be the case even if the feasible region is not compact). \\

\textbf{Initialization:} {Initial polyhedron $P^0$ has $m+1$ vertices, $V^0 = \{v^0, v^1, \ldots, v^m\}$, where $v^0 = y$, $v^i = v^0 + Mz^i$ for $i=1,\ldots,m$. Moreover, there are $l+1$ facets, $l$ of which correspond to the facets $\tilde{F}$ of cone $K$ in the sense that $f^i=(\{y\}+\tilde{f}^i) \cap P^0$, for $i=1,\ldots, l$, and $f^0$ is given by $\conv\{V^0\setminus\{v^0\}\}$. Then, $F_{v^i} =  \{f^0\}\cup \{f^j\st \tilde{f}^j\in \tilde{F}_{z^i}\}$ for $i=1,\ldots,m$ and $F_{v^0}=\{f^1,\ldots,f^l\}$. Moreover, $V_{f^0}=V^0\setminus\{v^0\}$ and $V_{f^i} = \{v^0\}\cup \{v_j \st z_j \in \tilde{V}_{\tilde{f}^i} \} $. }

\floatname{algorithm}{Variant}
\begin{algorithm}[!]
	\addtocounter{algorithm}{-3}
	\begin{algorithmic}
		\IF{$i=0$}
		\STATE Initialize $P^0$ as described (with vertices $V^0$ and facets $F^0$ and adjacency lists);
		\ELSE
		\STATE $P^i = onlinevert(P^{i-1},H)$;
		\ENDIF
	\end{algorithmic}
	\caption{Substitution of line 5 of Algorithm~\ref{alg:benson}}
	\label{alg:variant1_1} 
\end{algorithm}
\begin{algorithm}[!]
	\addtocounter{algorithm}{-1}
	\begin{algorithmic}
		\FORALL{$v\in V^i$}
		\IF{$v\in f^0$}
		\STATE $V^i \gets V^i \setminus \{v\}$;
		\ENDIF
		\ENDFOR
		\RETURN $\bar{\mathcal{X}}$: Set of weakly efficient solutions \\ ~~~~~~~~~~~~$P := \conv V^i + {K}$: Polyhedral outer approximation to the upper image.
	\end{algorithmic}
	\caption{Substitution of line 17 of Algorithm~\ref{alg:benson}}
	\label{alg:variant1_2} 
\end{algorithm}

The changes in the pseudo-code for this variant is given by Variant~\ref{alg:variant1_1}. Note that the vertices on facet $f^0$ are `artificial' by the construction of the initial polyhedron, hence the vertices on facet $f^0$ of $P^0$ are eliminated from the set of vertices of the current (last) polyhedron. This is why one needs additional lines before returning the output of the algorithm.  

\subsection{Variant 2: {VOP} Algorithm with Modified DD Method} \label{subsect:Alg2}
For this variant of Algorithm~\ref{alg:benson}, we call the modified DD method in order to solve the offline VEP for unbounded polyhedrons. The structure of the main algorithm is almost the same as the previous one. The only difference is in its initialization. The changes in the pseudo-code is given in Variant~\ref{alg:variant2}. \\

\textbf{Initialization:} Note that {$P^0 = \{y\} + K$ is the initial polyhedron. Then, the set of vertices of the initial polyhedron is $V^0=\{y\}$. Moreover, there are $l$ facets given by $f^i = \{y\} + \tilde{f}^i$, for $i=1,\ldots,l$. The adjacency lists are $F_{v^0} = \{f^1,\ldots,f^l\}$, $F_{z^i}= \{f^j \st \tilde{f}^j \in \tilde{F}_{z^i} \}$ and $V_{f^i}=\{v^0\}\cup \{z^j \st z^j\in \tilde{V}_{\tilde{f}^i} \}$ for $i=1,\ldots,l$. 
}
\begin{algorithm}[!]
	\begin{algorithmic}
		\IF{$i=0$}
		\STATE Initialize $P^0$ as described (with vertices $V^0$ and facets $F^0$ and adjacency lists);
		\ELSE
		\STATE $P^i = onlinevert2(P^{i-1},H)$;
		\ENDIF
	\end{algorithmic}
	\caption{Substitution of line 5 of Algorithm~\ref{alg:benson}}
	\label{alg:variant2} 
\end{algorithm}

\section{Computational Tests}\label{sec:numerical}

There is a MATLAB implementation of Algorithm~\ref{alg:benson}, which uses the vertex enumeration procedure (\emph{vert}) that has been used also in \cite{bensolve12}. The procedures explained in Section~\ref{sec:DD} as well as the variants of Algorithm~\ref{alg:benson} given in Section~\ref{sec:alg12} are implemented using MATLAB.

In order to compare the performances of the vertex enumeration procedures of Algorithm~\ref{alg:benson} and Variants~\ref{alg:variant1_1}-\ref{alg:variant2}, we randomly generate linear multiobjective optimization problems in the following form:
\begin{align*}
\text{minimize~~~} Cx \text{~~~subject to~~~} A x \leq b, \:\: x \geq 0, 
\end{align*}
where $C\in \R^{d\times n}$, $A \in \R^{m \times n}$, $b\in \R^m$ {and the ordering cone is $\R^d_+$}. As the objective space is $d$-dimensional, the vertex enumeration problem to be solved is also $d$-dimensional, see Algorithm~\ref{alg:benson} line 5. For our tests, each component of $A$ and $C$ is generated using independent normal distributions with mean $\mu=0$ and variance $\sigma^2 = 100$, whereas each component of vector $b$ is generated using independent uniform distributions on range $[0,10]$. In order to avoid numerical complications, we round each component of $A, C$ and $b$ to its closest integer. When we generate a problem, we first check the feasibility and boundedness (in the sense that the ideal point is finite) of it and add it to our sample only if the problem is bounded and feasible, hence solvable. Otherwise, we continue generating another set of $C, A$ and $b$.  

We generate different set of problems where the number of objectives ($d$) ranges from 2 to 4; the number of constraints ($m$) is taken as $2n$, where $n$ is the number of variables. For the problems with two objectives ($d=2$), we generate 30 feasible and bounded linear MOPs and for $d=3$ and $d=4$, we generate 20 of them. 

The aim of the computational tests is to compare the performances of different vertex enumeration procedures that is called in each iteration of Algorithm 1 (respectively Variant 1 and 2). Note that the efficiencies of Algorithm 1 and the two variants depend also on the choice of the vertex to be considered in each iteration, see line 7 of Algorithm 1. Indeed, for the current implementation of Algorithm 1, a vertex $v$ is chosen arbitrarily (depending on the structure of the list of vertices to be considered). Hence, calling Algorithm 1, and Variants 1 and 2 separately for the same test problem and checking the overall performances does not necessarily yield a fair comparison of the performances of the vertex enumeration procedures that is used within. It is possible that the algorithm and the variants go over the vertices of the current outer approximation in different orders. Hence, starting from the earliest iterations, each variant may yield a different current outer approximation, which of course affect the overall performance.   	

In order to overcome this difficulty, we solve the problems using Algorithm~\ref{alg:benson}, but in each iteration we solve the same vertex enumeration problem using three different methods: \emph{vert} from~\cite{bensolve12}, $onlinevert$ and $onlinevert2$. In order to do that we had three different initialization for each procedure. We set $M = 10^{4}$ for Variant 1. We measure the CPU times that is spent during each vertex enumeration procedures starting from the first iteration. The approximation error $\epsilon$ is taken as $0.005$ for the bi-objective problems and as $0.05$ for the problems with more than two objectives. 

The tests are conducted on a computer with system features Intel(R) Core(TM) i5-7200U CPU@ 2.50 GHz 2.71 GHz, 8.00 GB RAM, X64 Windows 10 and we utilize MATLAB R2013a.

\begin{figure}[h]
	\centering
	\includegraphics[width=7.4cm, height=4.3cm]{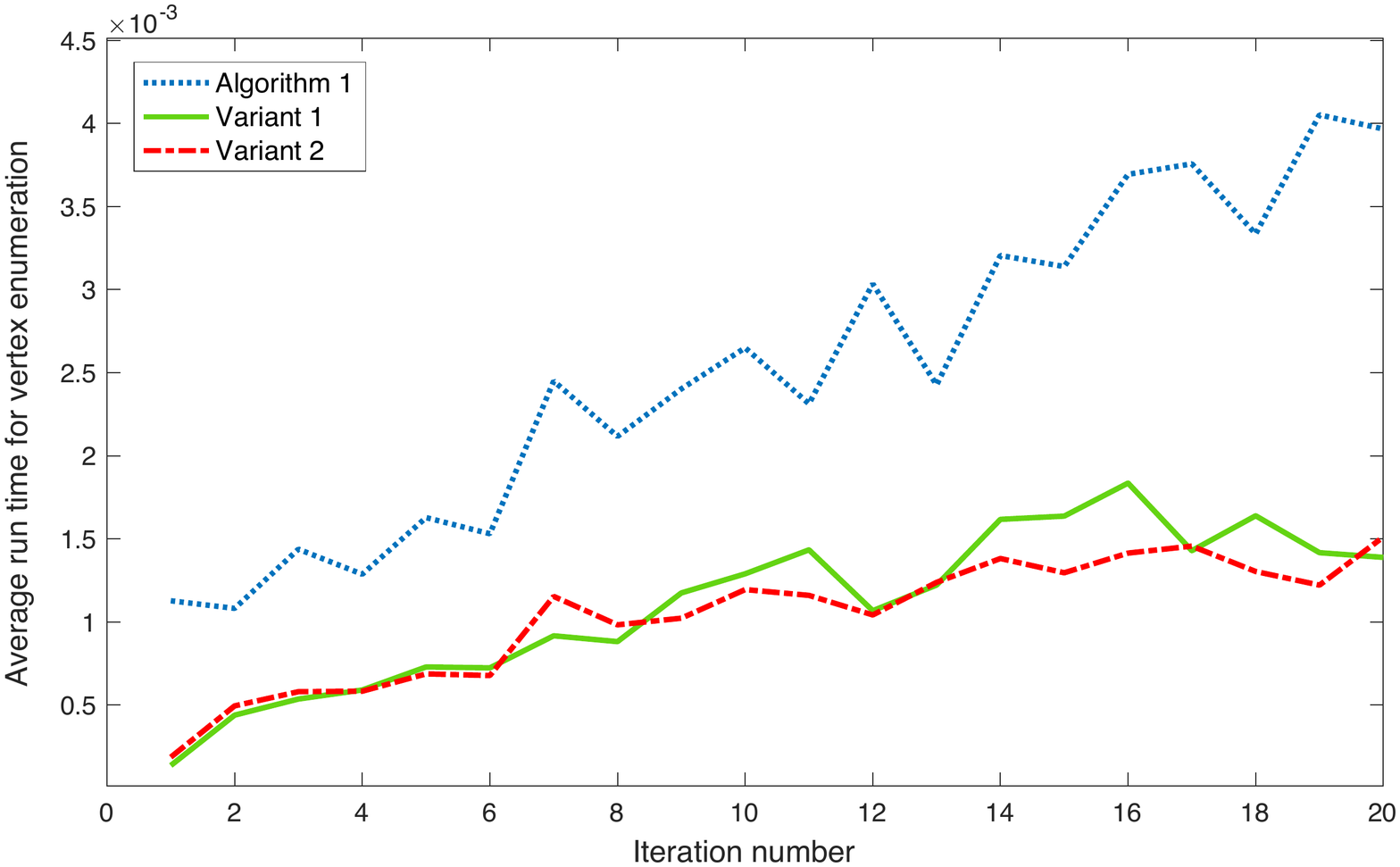} \\
	\caption{Run time performances of vertex enumeration procedures for problems with $d=2$ and $n=50$. The total number of iterations is 20; and the sample size is 20.}	
	\label{fig:2dim}
\end{figure}

We compare the average CPU times that are spent during these vertex enumeration procedures in each iteration of Algorithm 1. Indeed, we consider a sub-sample of problems: To explain it with an example, for $d=2$, $n = 50$, we generate $30$ problems among which the minimum number of iterations (of Algorithm 1) required is observed as $3$. If we want to have a sample of $30$ instances requiring the same number of iterations, we need to consider only the first 3 iterations of all these problems. Instead of considering $30$ rather small-sized (3 iterations) problems, we reduce the sample size to $20$ and we increase the number of iterations accordingly. More precisely, we list $30$ problems according to the number of iterations that they require in a non-increasing order. Then, we consider the first 20 problems. Figure~\ref{fig:2dim} shows the average run time spent for the vertex enumeration in each iteration. We observe that the run time of the vertex enumeration procedure used in Algorithm 1 is slightly more than the twice of the time spent by the Variants 1 and 2. However, there is no clear distinction between the two variants for these instances.

\begin{figure}[h]
	\centering
	\includegraphics[width=7.4cm, height=4.3cm]{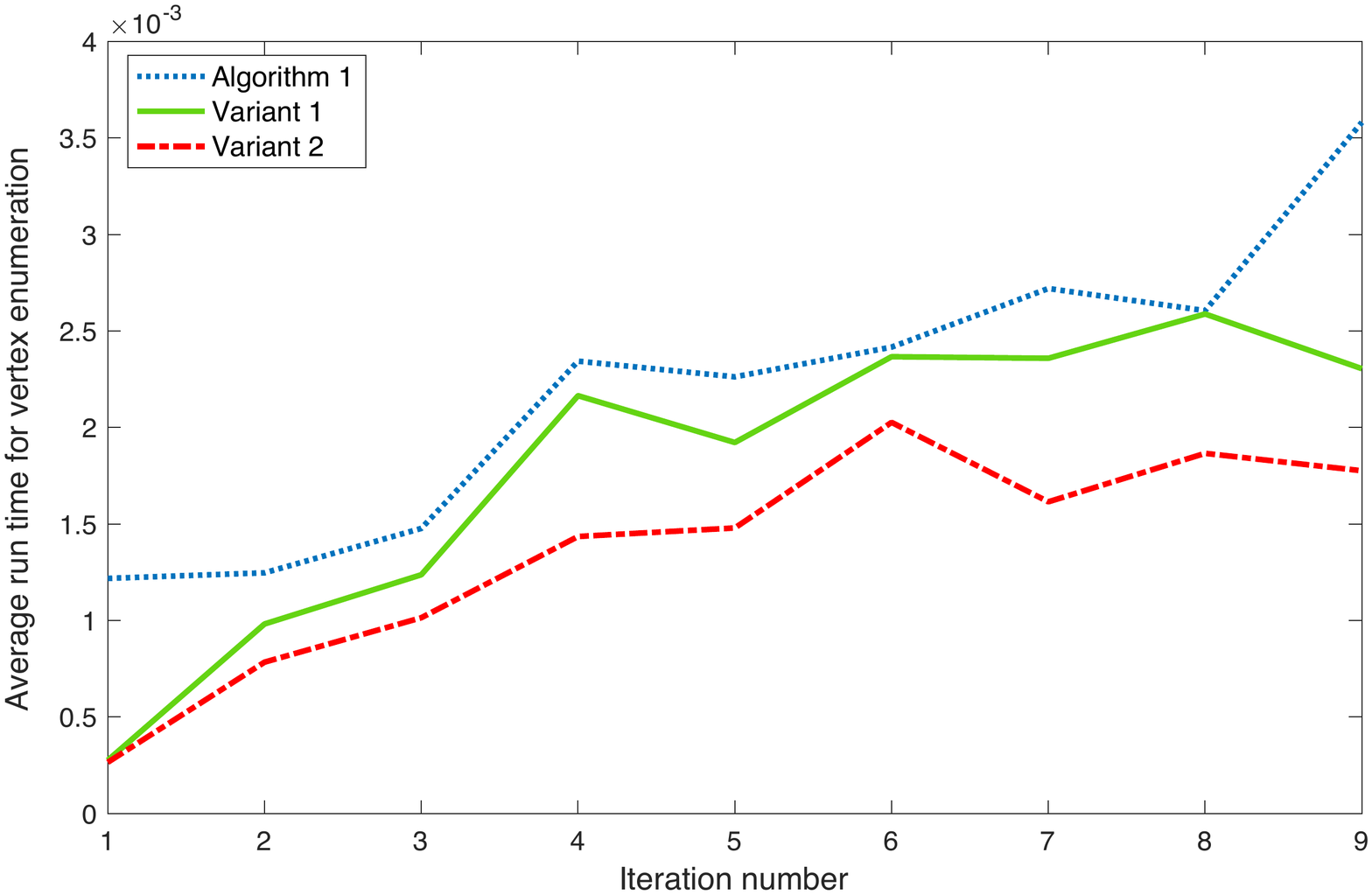}
	\includegraphics[width=7.4cm, height=4.3cm]{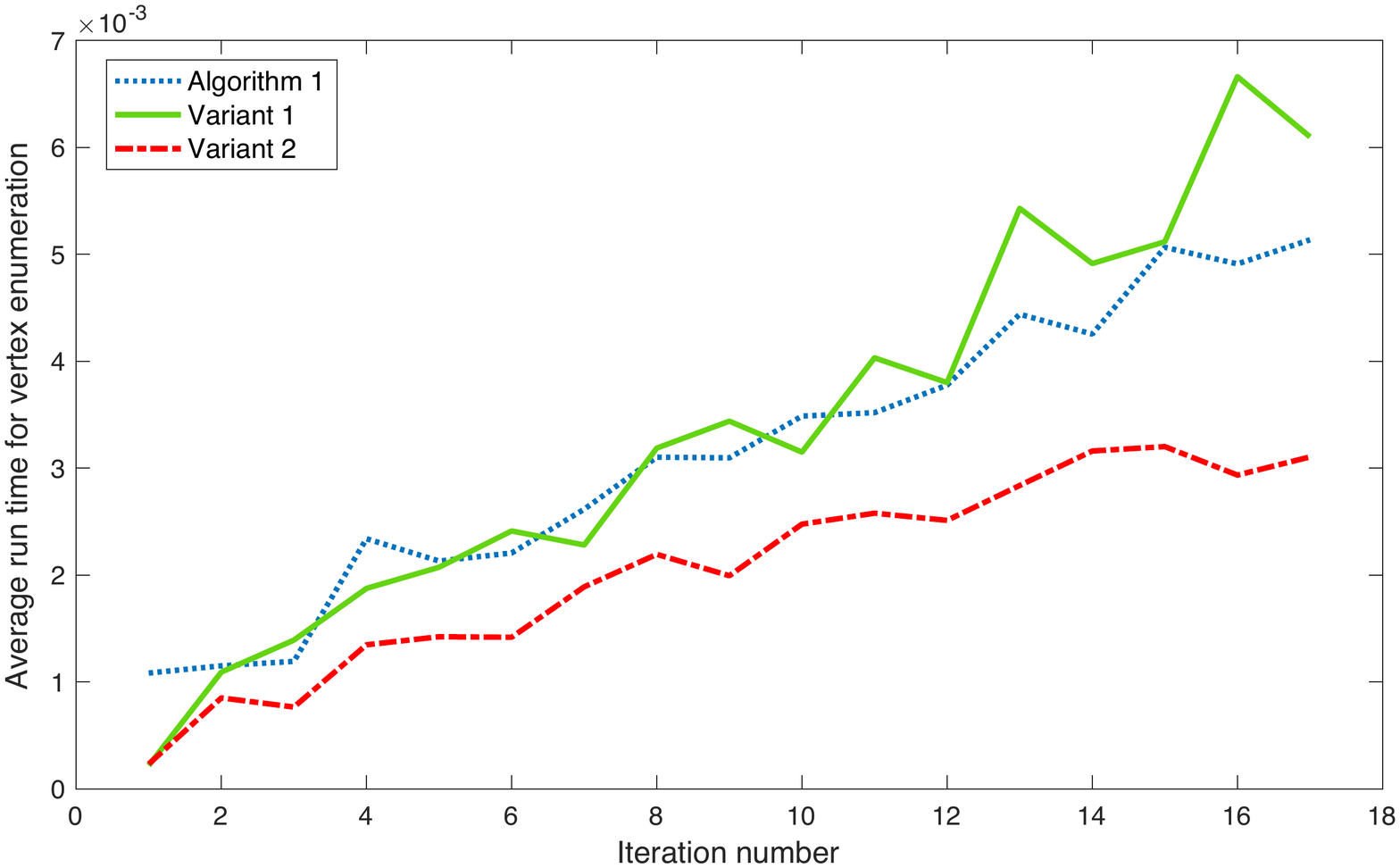}\\
	\includegraphics[width=7.4cm, height=4.3cm]{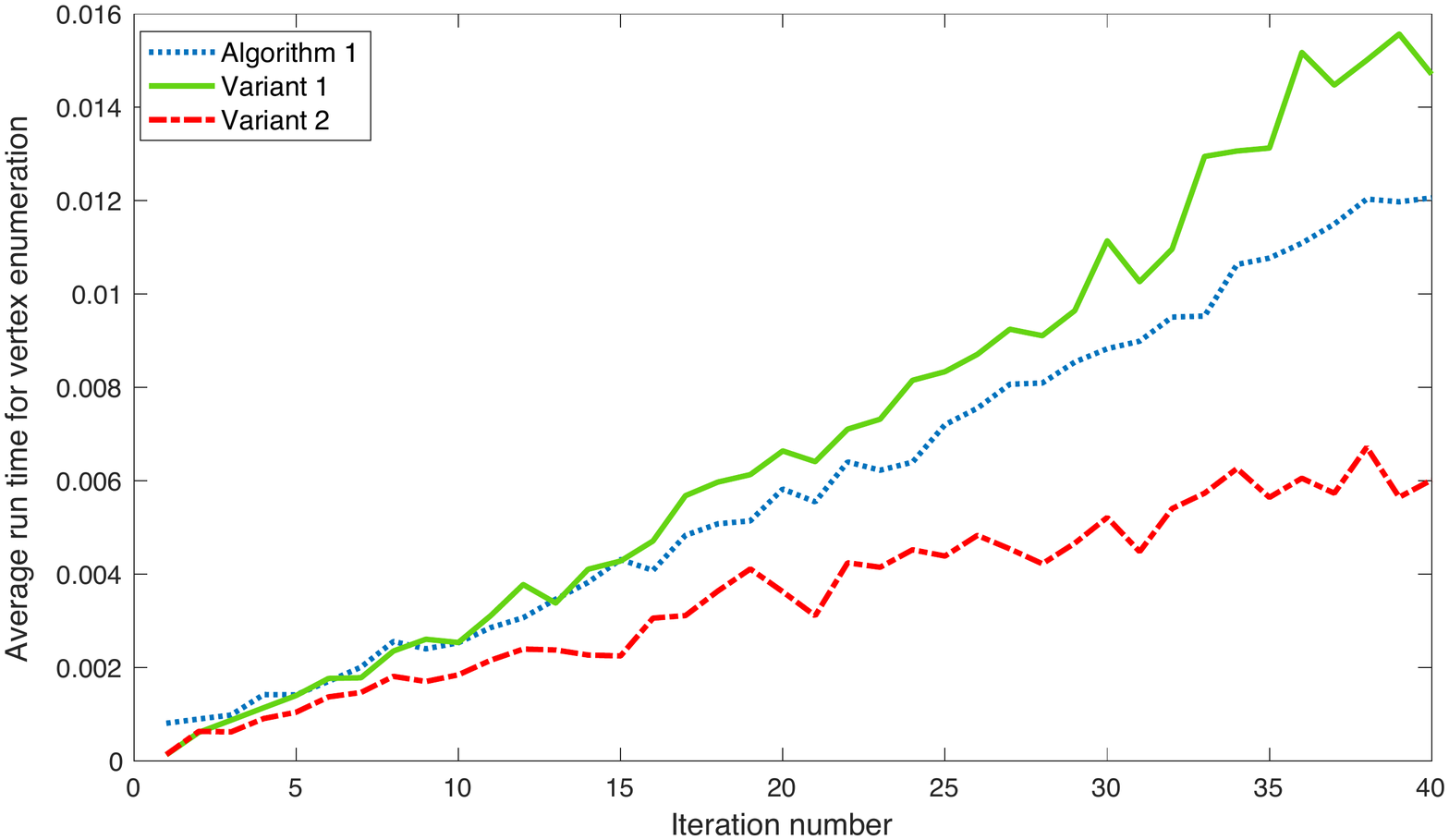}
	\includegraphics[width=7.4cm, height=4.3cm]{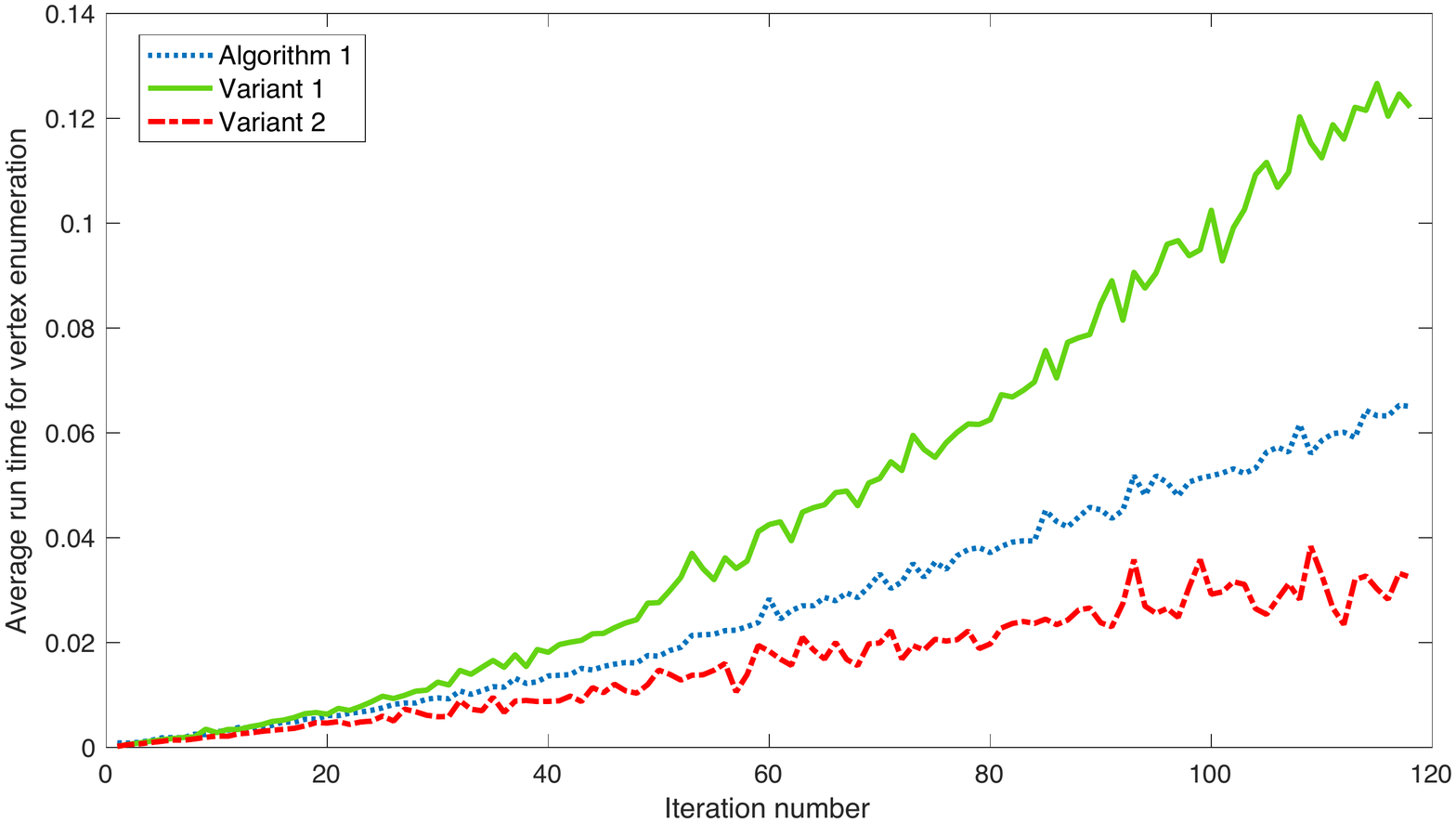}\\
	\caption{Run time performances of vertex enumeration procedures for problems with $d=3$ and $n=5$ (top left), $n=10$ (top right), $n=20$ (bottom left), $n=30$ (bottom right). The total number of iterations are 9, 17, 40, 118, respectively; and the sample size is 15 for all.}	
	\label{fig:3dim}
\end{figure}

For $d=3$, we consider four set of parameters, where we take the number of variables $n$ as 5, 10, 20, 30. Here, we expect that the MOPs would require more iterations as the size of the problem increases. The motivation of generating different sizes is to observe this pattern and if this is the case, then to observe the performance of the different vertex enumeration procedures as the iteration number increases. For each set, we generate $20$ problems and consider the sub-samples of sizes $15$, as explained before. The graphs can be seen in Figure~\ref{fig:3dim}. 

As the number of variables of the MOP increases, we observe that the number of iterations required for the algorithm increases, as expected. Moreover, we see that the average CPU time spent for each iteration increases as the iteration number increases. This is expected since, in general, the number of vertices to be checked in each iteration increases. 

In the first two graphs ($n=5$ and $n=10$) of Figure~\ref{fig:3dim}, the performances of Algorithm 1 and Variant 1 are similar, whereas Variant 2 seems to work faster then both. As the iteration number increases, which is the case for larger problems ($n=20$ and $n=30$), the differences in the performances also increase. Moreover, it is observed that Variant 1 gets worse than Algorithm 1 as the iteration number increases. 

\begin{figure}[!]
	\centering
	\includegraphics[width=7.4cm, height=4.3cm]{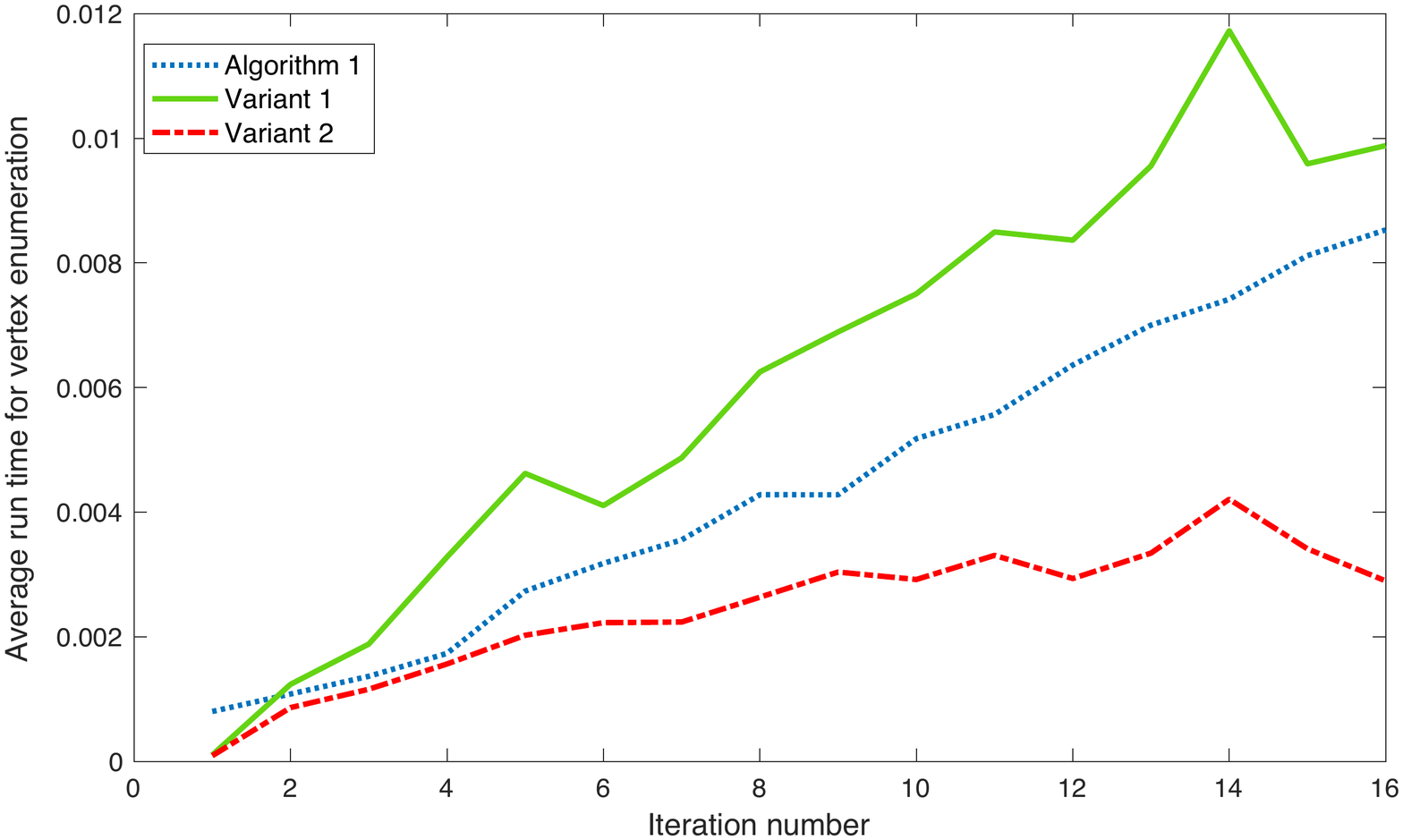} 		
	\includegraphics[width=7.4cm, height=4.3cm]{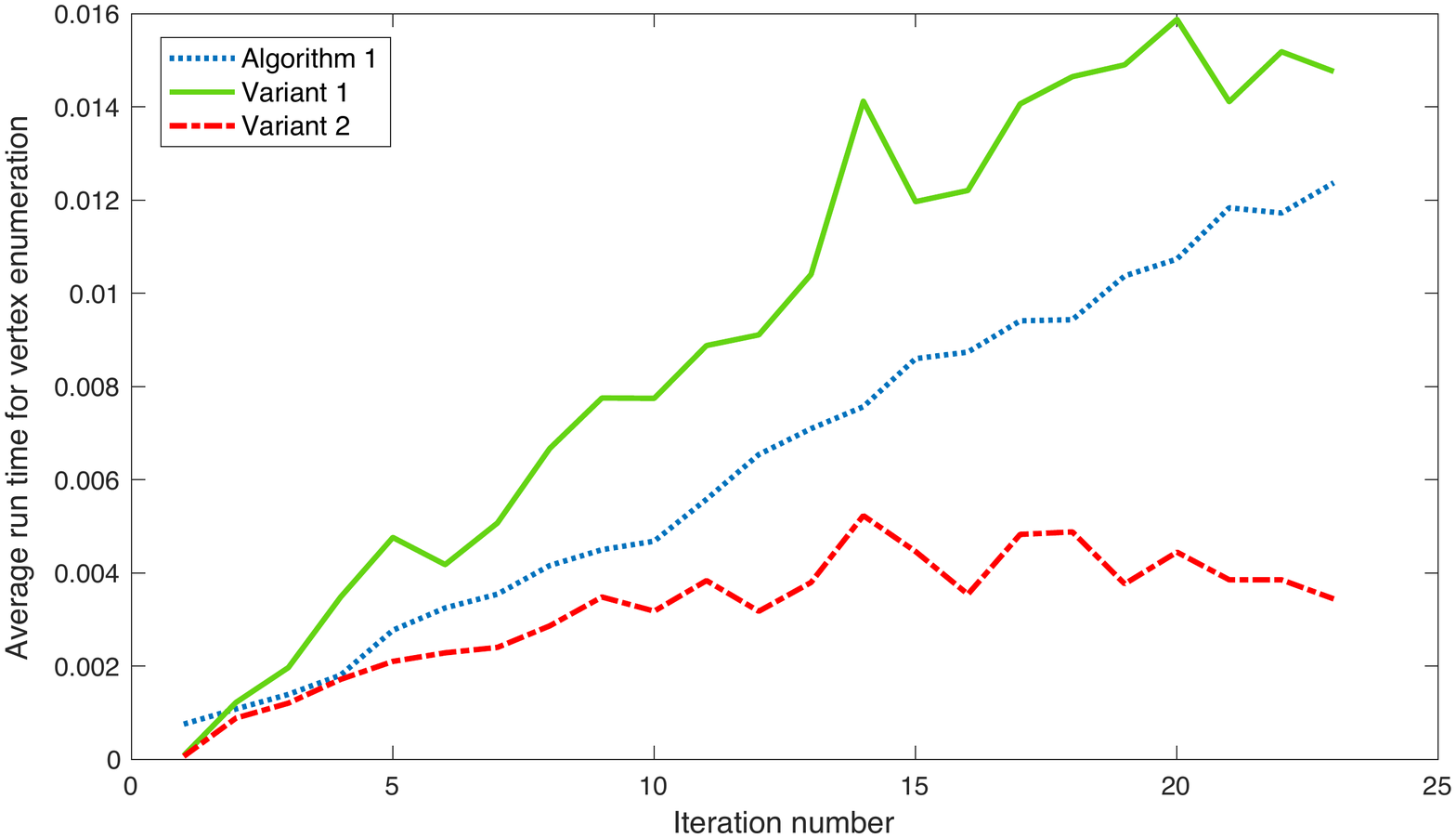}		
	\caption{Run time performances of vertex enumeration procedures for problems with $d=4$, $n=5$. The total number of iterations are 16 and 23; and the sample size is 15 and 10, respectively (left and right).}	
	\label{fig:4dim}
\end{figure}

For $d=4$, we generate 20 problems with $n=5$ and we consider two sub-samples of sizes $15$ and $10$. The graphs can be seen in Figure~\ref{fig:4dim}. We observe a similar pattern as we observed for three dimensional problems. Different from those, the pattern is clear even from the earliest iterations.

By checking Figures~\ref{fig:2dim}-\ref{fig:4dim}, we observe that the time requred for vertex enumeration increases as the dimension of the (objective) space $d$ increases. For example, if one considers the average CPU time spent in $20^{th}$ iteration for $d = 2, d = 3$ (check for instance bottom left figure with $n=20$) and $d=4$, these are respectively around $0.004, 0.006$ and $0.01$ for Algorithm 1; $0.0015, 0.007$ and $0.015$ for Variant 1; and $0.0015, 0.003$ and $0.004$ for Variant 2. Indeed, we see that the increase in the run times is the most for Variant 1. 

Note that \emph{vert} from~\cite{bensolve12} solves an off-line VEP whereas \emph{onlinevert} solves an on-line VEP in each iteration. Hence, it may not be expected to observe that \emph{vert} (Algorithm 1) is more efficient than \emph{onlinevert} (Variant 1), especially as the number of iterations increases. This occurred for higher dimensional problems possibly because there are too many artificial vertices to be considered even though they are deleted at the very end of the algorithm, see Variant~\ref{alg:variant1_2} - Substitution of line 17 of Algorithm~\ref{alg:benson}. For two sets of random instances (15 instances with $d=3, n=20$ and 10 instances with $d=4, n=5$), we check the number of artificial and actual vertices that are generated at each iteration. Among these random instances, the minimum numbers of iterations required are 45 for the set with $d=3$; and 17 for the set with $d=4$. Hence, we check the first 45 and respectively, 17 iterations of corresponding sets of instances. Each row in Table~\ref{tab:artif} shows for the particular iteration, the average number of actual and artificial vertices as well as the percentage of the artificial vertices within all. Note that instead of providing this information for every iteration, we pick some of them as this is sufficient to summarize the general trend. 

For the problems with $d=3$, we observe that on the average, 20 percent of the vertices considered for \emph{onlinevert} in each iteration are artificial. This percentage is higher for the earlier iterations and decreases later on. However, the average number of artificial vertices are increasing as in general the number of vertices increases rapidly through iterations. Indeed, for a sub-sample of size 10, we can increase the iteration number up to 70 and we observe the same pattern, see the last three rows of Table~\ref{tab:artif}.

For the problems with $d=4$, more than half of the vertices considered for $onlinevert$ are observed to be artificial. Different from the previous case, this percentage is increasing through the iterations. Note that we can increase the iteration number up to 26 by decreasing the sample size (see the last three rows of Table~\ref{tab:artif}) and the same pattern holds even then. This explains the poor performance of Variant 1 for high dimensional problems.  

\begin{table}[h] 
	\centering
	\resizebox{\textwidth}{!}{
		\begin{tabular}{|c|c|c|c||c|c|c|c|}
			\hline
			\multicolumn{4}{|c||}{$d=3, n=20$} & \multicolumn{4}{c|}{$d=4, n=5$} \\ \hline
			\# iteration & \begin{tabular}[c]{@{}l@{}} \# actual\\ vertices\end{tabular} & \begin{tabular}[c]{@{}l@{}} \# artificial\\  vertices\end{tabular} & \% artificial & 		\# iteration & \begin{tabular}[c]{@{}l@{}} \# actual\\ vertices\end{tabular} & \begin{tabular}[c]{@{}l@{}} \# artificial\\  vertices\end{tabular} & \% artificial \\  \hline	
			10  & 14.00 & 7.27   & 34.17 & 5   & 4.80  & 7.30  & 60.33  \\
			20 & 27.27 & 10.80 & 28.37 & 10 & 12.10  & 17.60 & 59.26  \\
			30 & 41.60 & 13.20 & 24.09 & 13 & 14.40	& 21.00 & 59.32 \\
			40 & 54.13 & 16.47 & 23.32 & 17 & 15.70	 & 25.20 &61.61 \\
			50$^*$ & 73.00  & 18.10	& 19.87 & 21$^{**}$ & 17.22	& 28.89 & 62.65 \\
			60$^*$ & 86.70  & 20.90 & 19.42 & 23$^{**}$ & 16.87	& 29.63 &	63,71 \\		
			70$^*$ & 100.90 & 24.30	& 19.41 & 26$^{**}$ & 17.00	& 32.50 & 	65.66\\ \hline
	\end{tabular}}
	\caption{Number of actual and artificial vertices for test instances with $d=3$, $n=20$ and with $d=4, n=5$. The sample size is 15 for $d=3$ instances (* except for the last three rows, for which the sample size is 10); and it is 10 for $d=4$ instances (** except the last three rows, for which the sample sizes are 9, 8 and 6, respectively).}
	\label{tab:artif}
\end{table}

\section{Conclusion}\label{sec:concl}

We study the vertex enumeration problem, in particular the DD method (\emph{onlinevert}) to be used within an objective space based {VOP} algorithm (Algorithm 1), which currently employs an offline vertex enumeration procedure \emph{vert}. We {also consider} a modified DD method (\emph{onlinevert2}) which works for unbounded polyhedrons {with known recession cones, as proposed in \cite{burton2010projective}. We provide a detailed Euclidean space description of this method}. We examine two variants of Algorithm 1, one using $onlinevert$ and the other using $onlinevert2$. We compare the performances through randomly generated linear MOP instances.

Overall, \emph{onlinevert2} used in Variant 2 is observed to be the most efficient procedure among the others especially as the dimension of the objective space and also as the number of iterations increases. {We demonstrate through a detailed computational study that}, for vertex enumeration problems, where the polyhedron to be computed is unbounded {with a known recession cone}, employing the {modified} variant of the DD method (\emph{onlinevert2}) {increases} the overall efficiency {significantly}. As discussed throughout, one crucial application area is the objective space based MOP algorithms. However, it could be employed as a subroutine for any algorithm or procedure that requires solving VEPs with the aforementioned property.

	\bibliographystyle{plain}
\bibliography{whole_paper}

\end{document}